 \documentstyle [12pt,twoside, epsf]{article}
\newcommand{\ZZ}{\mbox{$Z\!\!\! Z\!$}}	
\newcommand{\QQ}{\mbox{$Q\!\!\!\! I$}}	
\newenvironment{pf}{\noindent {\em Proof.}}{\hfill $\Box$ \medskip}

\newtheorem{thm}{Theorem}[section]
\newtheorem{lem}[thm]{Lemma}
\newtheorem{cor}[thm]{Corollary}
\newtheorem{prop}[thm]{Proposition}

\newtheorem{defn}[thm]{Definition}

\newtheorem{exmp}[thm]{Examples}

\newtheorem{rems}[thm]{Remarks}

\def\picill#1by#2(#3)
{\vbox to #2
{\hrule width #1 height 0pt depth 0pt
\vfill\epsffile{#3}}}

\let \ttorg \tt \def \tt{\ttorg \obeyspaces}

\begin{document}
\pagestyle{myheadings}
\markboth{Geck--Lambropoulou}{Markov traces and knot invariants}

\title{Markov traces and knot invariants related to Iwahori-Hecke 
algebras of type $B$} 
\author{Meinolf Geck and Sofia Lambropoulou} 

\date{}

\maketitle


\section{Introductory remarks on knots, braids and trace functions}

In  classical knot theory we  study knots inside the $3$-sphere modulo
isotopy.   Using the Alexander and Markov  theorem,  we can translate
this into a purely  algebraic setting in terms  of Artin  braid groups
modulo an equivalence  relation generated  by  `Markov moves'  (one of
which   is usual  conjugation  inside  the   braid
group).    V.F.R.\  Jones \cite{Jones} used  this    fact  in 1984 for
constructing   a new knot   invariant  through trace functions  on the
associated Iwahori-Hecke  algebras of type  $A$  with suitable 
properties that reflect  the  above  Markov moves.  

Jones's work led to  questions of developing knot theory corresponding
to other types of Coxeter groups. It is proved
in \cite{L1}  that there exist braid   structures related to arbitrary
3-manifolds, which   in  addition  satisfy  appropriate Markov-isotopy
equivalence; also that, if the 3-manifold is a solid torus,
then  the sets of  related braids form groups, which   are in fact the
Artin-Tits braid groups related  to the $B$-type Coxeter groups. These
results together with a linear  trace that we  found in 1991 are
used  in  \cite{L2} for constructing    a  4-variable analogue of  the
homfly-pt (2-variable Jones)   polynomial for oriented knots  inside a
solid torus.  We proved the existence of this trace (see \cite{L2}) by
following and adapting   to the $B$-type  case  Jones's proof of   the
existence of Ocneanu's trace in \cite{Jones}, Theorem~5.1.

\bigbreak
The aim of this  paper is to give a  full classification of  {\it all}
linear traces on Iwahori-Hecke algebras  of type $B$ which support the
{\it Markov property}, see  Definition~4.1 and Theorem~4.3.  This uses
in an essential way the results in  \cite{GP} about trace functions on
arbitrary  Iwahori-Hecke  algebras associated  with   finite   Coxeter
groups.  This method yields an  alternative proof of  the above
special trace (cf.\ also  \cite{GP}, (4.2), where an alternative proof
for Ocneanu's original trace is given.) 

In Section~5 we  discuss  the knot  theory of a   solid torus  and  we
explain why the related braid groups are in fact the Artin-Tits groups
of $B$-type.  We also give the  Markov equivalence of these braids, so
that the equivalence classes correspond bijectively to isotopy classes
of knots  in  the solid torus (detailed  account  and proofs of  these
results     can be  found  in   \cite{L1}  or  \cite{L2}). The  Markov
equivalence is in   terms of two   isotopy moves  which are  reflected
precisely in  the definition of the  Markov  property for  our traces.
Then we normalize  properly the constructed traces  in order to obtain
{\it all} homfly-pt analogues   related to the  $B$-type Iwahori-Hecke
algebras,    for    oriented    knots  inside  the   solid
torus.  Finally, we give the  skein rules  and initial conditions that
characterize these invariants diagramatically.

In \cite{G} the first author uses the results of this paper to provide a full
classification of Markov traces for  Iwahori-Hecke algebras of  type $D$ (see 4.7
for precise statement of the main result).  Moreover, in a  further  `vertical' 
development  (cf. \cite{L3}) the second author considered {\it all} Hecke-type
quotients of  the Artin-Tits  braid   group of $B$-type  and  constructed Markov
traces and knot  invariants on all levels, the  basic  level being the Iwahori-Hecke
algebras of $B$-type and the results in \cite{L2}.

\bigbreak We shall now explain  in more detail  our results.  Let us
consider the following Dynkin diagram.

\begin{picture}(300,50)
\put(  0, 20){$(B_{n})$}
\put( 40, 20){\circle*{5}}
\put( 40, 18){\line(1,0){20}}
\put( 40, 22){\line(1,0){20}}
\put( 60, 20){\circle*{5}}
\put( 60, 20){\line(1,0){30}}
\put( 90, 20){\circle*{5}}
\put( 90, 20){\line(1,0){20}}
\put(120, 20){\circle*{1}}
\put(130, 20){\circle*{1}}
\put(140, 20){\circle*{1}}
\put(150, 20){\line(1,0){20}}
\put(170, 20){\circle*{5}}
\put( 37, 30){$t$}
\put( 57, 30){$s_1$}
\put( 87, 30){$s_2$}
\put(167, 30){$s_{n-1}$}
\put(250, 20){$n \geq 1$}
\end{picture}

\noindent The symbols $t,s_1,\ldots,s_{n-1}$ labelling the nodes form
generators for  the corresponding Artin-Tits braid group $\tilde{W}_n$
and the finite Coxeter group $W_n$ of type $B_n$. The braid group 
$\tilde{W}_n$ has the defining relations 
\[\begin{array}{rclll}
s_1ts_1t & = & ts_1ts_1  & & \\
ts_i & = & s_it & \mbox{ if } & i>1 \\
s_is_j & = & s_js_i & \mbox{ if } & |i-j|>1 \\
s_is_{i+1}s_i & = & s_{i+1}s_is_{i+1} & \mbox{ if } & 1 \leq i \leq n-2
\end{array}\]
Relations of these  types will  be called {\it  braid relations}.   As
proved in \cite{L1} the elements  of $\tilde{W}_n$ can be  represented
geometrically by braids in  $S^3$ on $n+1$ strands  in which the first
strand remains  pointwise  fixed.  When  we  refer to  this  geometric
interpretation of $\tilde{W}_n$ we shall denote it by $B_{1,n}$. Below
we illustrate the generators $s_i, t$ and the element $t_i'=s_i \cdots
s_1ts_1^{-1} \cdots s_i^{-1}$ in $B_{1,n}$, and  also an example of an
element in $B_{1,5}$.

 $\vbox{\picill8inby1.75in(mein1)  }$

\noindent If in addition to the braid relations we impose the quadratic
relations that  each generator has order~$2$, then  we obtain the finite factor
group  $W_n$. 
\bigbreak
The corresponding Iwahori-Hecke
algebra $H_n$ is obtained as a   quotient of the group  algebra   of $\tilde{W}_n$
by factoring out the quadratic  relations  
\[ t^2=(Q-1)t+Q \cdot 1  \quad \mbox{ and } \quad g_i^2=(q-1)g_i+q \cdot 1 
\mbox{ for all $i$},\]  
where we denote the image of $s_i$ in  $H_n$ simply by $g_i$, and where $q, Q$ are
fixed parameters from the ground ring. The  algebra $H_n$ is finite-dimensional,
with a basis $\{g_w\}$ labelled by the elements of $W_n$. Now the idea  is to
construct invariants  of knots in the solid torus using trace    functions on
$\tilde{W}_n$  which  factor through $H_n$ and which   respect  the braid
equivalence on  $\tilde{W}_n$. The latter is generated by the following two moves 
(cf.\ Theorem~5.2).  \begin{itemize}
\item[(i)] Conjugation: if $\alpha,\beta \in \tilde{W}_n$
then $\alpha \sim \beta ^{-1}\alpha\beta$.
\item[(ii)] Markov moves: if $\alpha \in \tilde{W}_n$ then $\alpha \sim
\alpha{s_n}^{\pm 1}\in \tilde{W}_{n+1}$. 
\end{itemize} 
\noindent Thus, the problem is  reduced to studying trace
functions $\tau$ on $H:= \bigcup_{n=1}^{\infty} H_n$ which  satisfy the rule 
$\tau(hg_n)=z\tau(h)$, where $z$ is a fixed parameter  in the ground  ring over
which the algebra $H$ is defined, and $h \in H_n$. (Note that such an $h$ is  a
linear combination of basis elements $g_w$ which do  not involve the  generator
$g_n$.) This rule is what  we call `the Markov property' for trace functions on
$H$.
\bigbreak
A general  scheme for constructing  trace functions on $H_n$ (in fact,
for  Iwahori-Hecke algebras of any given  type)  has been developed in
\cite{GP}. Firstly, it is known that any trace is determined by its values on basis
elements corresponding to a set of representatives of the conjugacy classes of
$W_n$. However, it is not true in general that basis elements
 corresponding to conjugate group elements are also conjugate in the algebra, and so,
to compute the trace of an arbitrary element is no more a trivial task.  In
\cite{GP}, there is an explicit algorithm  for computing the value of a given trace
on an arbitrary basis element $g_w$ from the values on basis elements corresponding
to elements  of {\it minimal length} in the various conjugacy classes. 
\bigbreak
\noindent In $W_n$, let
$t_i:=s_i \cdots s_1t s_1 \cdots s_i$ and let a
positive, respectively  negative,  block be an element of the form
\[s_{i+1}s_{i+2}\cdots s_{i+m} \quad \mbox{ respectively } \quad
t_i s_{i+1} \cdots s_{i+m}, \quad \mbox{ where } m,i\geq 0.\]
Then an element of minimal length in a conjugacy class is a  product of negative
blocks ordered by increasing length, followed by various positive blocks also in 
increasing  length (for example, $t t_1s_2 t_3s_4 s_6 s_8s_9s_{10} \in W_{11}$). The
main idea of this paper is that, in order to determine the  
 Markov traces on  $H$, we lift the elements $t_i$ in $W_n$ to  the elements
$t'_i:=g_i \cdots g_1t g_1^{-1} \cdots g_i^{-1}$ in $H_n$ (instead of the obvious
lifting to $g_i \cdots g_1t g_1 \cdots g_i$). This will enable us to
parametrize in Section 4 all possible Markov traces with parameter~$z$  by the
initial conditions  \[ \tau(t_0't_1't_2' \cdots t_{k-1}')=y_k \quad \mbox{ for all
$k \geq 1$,}\] where $y_1,y_2 \ldots$ are arbitrary elements from the ground ring. In
particular, if $h=d_1 \cdots d_n \in H$ with $d_i \in \{1, g_{i-1}, t_{i-1}'
\}$, a lifting of a minimal length representative, we define a  trace function
$\tau$ on $H$ by the rule $ \tau(h)=z^{a(h)} y_{b(h)},$ where $a(h)$  is the 
 number of $d_i$ which are  in  the set $\{g_1,\ldots,g_{n-1}\}$ and $b(h)$ is the
number  of $d_i$ which are of the form $t_i'$ (cf. Theorem 4.3). From our definition 
it is clear that $\tau$ satisfies  the Markov property for  these
elements,  and the  problem is then to show that this  property holds on
all elements of $H$. This will require some technical preliminaries which  
are provided in Sections~2 and~3. The final proof will then be given in Section~4.

\section{Computations in the braid group of type $B_n$}

It is the purpose of  this section to reformulate  some of the results
in \cite{GP}  on elements in signed block  form in  terms of the braid
group. These will then carry through to the Iwahori-Hecke algebra level and will be
used in the  existence and  uniqueness proof for the  analogues of   Ocneanu's 
trace for type $B$.
 
\medskip
2.1. Let $W_n$ and $\tilde{W}_n$ as in Section 1. For  each element $w$ in $W_n$ (or
in $\tilde{W}_n$) we  define the length, $l(w)$, to be the smallest non-negative 
integer $k$ such that $w$  can   be  written as  a    product of $k$  generators  
(or their inverses).  Such an expression of  $w$ of minimal possible length will
be called  a reduced expression for  $w$. (See \cite{Bour}, Chap.\ IV,
\S 1.1.) The  exchange  condition for Coxeter  groups implies   that, if we are
given two  reduced expressions of an element  in $W_n$ as  products in
the generators    $t,s_1,\ldots,s_{n-1}$      then the   corresponding
expressions  in the braid  group   $\tilde{W}_n$ are also  equal  (see
\cite{Bour}, Chap.\ IV, \S 1, Proposition~5).

 By convention, we let $W_0=\{1\}$. Then, for  all $n \geq 1$, the
group $W_{n-1}$ is a parabolic subgroup  of $W_n$ obtained by removing
the  node  with   label $s_{n-1}$   (cf.\ \cite{Bour}, Chap.\   IV, \S
1.8).  Hence, we have natural  embeddings $W_0 \subset W_1 \subset W_2
\subset \ldots$ and we let $W:=\bigcup_{n} W_n$.  We define
\[t_i:=s_is_{i-1} \cdots s_1ts_1 \cdots s_{i-1}s_i \in W_n \quad
\mbox{ for all $0 \leq i \leq n-1$, where $t_0:=t$}.\] Then the set of distinguished
right  coset representatives of $W_{n-1}$ in $W_n$ is given as follows.
\[ \begin{array}{lcl}
{\cal R}_n:=\{ & 1,\;\;\; t_{n-1}, & \\
         & s_{n-1}s_{n-2} \cdots s_{n-k} & (1 \leq k \leq n-1), \\
         & s_{n-1}s_{n-2} \cdots s_{n-k}t_{n-k-1} & (1 \leq k \leq n-1)\}\\
\end{array}\]
(Note that ${\cal R}_1=\{1,t\}$.) Then each element $w \in W_n$ can be
written uniquely in  the form $w=r_1 \cdots  r_n$ with  $r_i \in {\cal
R}_i$.   Such an   expression  of $w$  is  reduced, that  is,  we have
$l(w)=l(r_1)+\ldots +l(r_n)$.

Finally,  let   ${\cal  D}_n:=\{1,s_{n-1},  t_{n-1}\}  \subseteq {\cal
R}_n$.  Then ${\cal  D}_n$  is the set  of distinguished  double coset
representatives of $W_{n-1}$ in $W_n$ (see \cite{Bour}, Chap.\ IV, \S 1, 
Ex.\ 3).

Each  $r  \in {\cal R}_n$  can  now be written  uniquely   in the form
$r=dr'$ where   $d \in {\cal  D}_n$  and  $r'=1$  or $r'=s_{n-2}\cdots
s_{n-k}$ or $r'=s_{n-2} \cdots s_{n-k}t_{n-k-1}$. In particular, we
have $r' \in {\cal R}_{n-1}$.

\medskip
2.2.  We   shall  now lift these elements to the   braid  group
$\tilde{W}_n$.  First  note   that we also     have natural embeddings
$\tilde{W}_0  \subset \tilde{W}_1 \subset  \tilde{W}_2 \subset \ldots$
and we  let $\tilde{W}:=\bigcup_{n}  \tilde{W}_n$.   Geometrically the
embedding of $B_{1,n}$ into $B_{1,n+1}$  is described by the following
picture.

$$\vbox{\picill3.7inby1.04in(mein2)  }$$   

\noindent We would like to define the analogue of $t_i$ to be a conjugate  
of $t$ where the conjugating element is of the form $s_1^{\pm 1} \cdots
s_i^{\pm 1}$. In accordance with the geometric considerations in Section~1,
we choose the exponents $\pm 1$ so as to obtain the element $t_i'$ 
already encountered above:
\[ t_i':=s_is_{i-1} \cdots s_1ts_1^{-1} \cdots s_{i-1}^{-1}s_i^{-1} 
\in \tilde{W}_n \quad \mbox{ for all $0 \leq i \leq n-1$}.\] 
(All inverses on the right hand side of $t$.) Then each $t_i'$ maps to
$t_i$  under  the canonical surjection $\tilde{W}_n  \rightarrow W_n$.
We  let ${\cal D}_n' \subseteq {\cal  R}_n'  \subseteq \tilde{W}_n$ be
the analogous sets as above, where each $t_i$ is replaced by $t_i'$.

For  any $i,j$ the  elements $t_j$  and $t_i$  commute with each other
(cf.\ \cite{GP}, (2.3)). For $t_j'$ and $t_i'$, this will only be true 
up to possibly changing some inverses in the definition of these elements. 
We will not completely formalize this, but only give the following 
additional definition.  For any $j \in \{0,\ldots,i\}$ we denote 
\begin{eqnarray*}
t_{i,j}' & :=  & s_i^{\pm 1} \cdots  s_{j+1}^{\pm 1} s_js_{j-1} \cdots
s_1ts_1^{-1} \cdots    s_{j-1}^{-1}s_j^{-1}s_{j+1}^{\mp      1} \cdots
s_i^{\mp 1}\\ &     =     & s_i^{\pm  1}   \cdots    s_{j+1}^{\pm   1}
t_j's_{j+1}^{\mp 1}  \cdots  s_i^{\mp  1}.\end{eqnarray*}  (That is,  all
inverses up   to index $j$  are in  the right  position,  and for each
bigger  index, the inverse may  be put either on  the  right or on the
left hand side.)

Finally,   we  let ${\cal  D}_{n,j}'$  be  the set  consisting of $1$,
$s_{n-1}$ and  all possible elements of  the form $t_{n-1,j}'$. 
Similarly, we define ${\cal R}_{n,j}'$.  As a convention, we will 
usually denote the elements in ${\cal D}_{n,0}'$ by the symbol $d_n^{*}$. 

The next result will show  that, in particular,  a product of the form
$t_i't_j'$ with $i < j$ can be written as $t_{j,i}'t_j'$. Relations of
this kind will be used frequently in the sequel.

\addtocounter{thm}{2}
\begin{lem}
The following relations hold in $\tilde{W}_n$.

\begin{itemize}

\item[(a)] $s_it_m'=t_m's_i$ and $s_i^{-1}t_m'=t_m's_i^{-1}$ for all 
$i<m$ and $i>m+1$.

\item[(b)] $t_i't_m'=s_m \cdots s_{i+2}s_{i+1}^{-1}s_i\cdots s_1ts_1^{-1}
\cdots s_i^{-1}s_{i+1}s_{i+2}^{-1}\cdots s_{m}^{-1}t_i'$ for $i < m$.\\
(The inverse at index $i+1$ changes in $t_m'$) 
\item[(c)] $s_{m-1} \cdots s_{m-k}t_{m-k-1}'t_{m}'=
t_{m,m-1}'s_{m-1} \cdots s_{m-k}t_{m-k-1}'$ for $0 \leq k \leq m-1$.
(The inverse at index $m$ changes in $t_m'$)
\end{itemize}
\end{lem}

\begin{pf} The defining relations for $\tilde{W}_n$ imply that
\begin{eqnarray*}
ts_i^{-1}=s_i^{-1}t & \mbox{  if } &  i>1 \\ s_is_j^{-1}=s_j^{-1}s_i &
\mbox{ if  } &  |i-j|>1 \\ s_is_{i+1}s_i^{-1}=s_{i+1}^{-1}s_is_{i+1} &
\mbox{ if  } & 1 \leq  i \leq n-2  \\ s_1^{-1}ts_1t=ts_1ts_1^{-1}. & &
\end{eqnarray*} The assertions   of  the Lemma now  readily  follow by
straightforward computations.   (Notice that these relations  could be
alternatively checked easily using the geometric interpretations given
above.) \end{pf}

\medskip
2.4. Consider an element of the form  $d_1 \cdots d_n \in \tilde{W}_n$
with   $d_i \in {\cal   D}_i'$ for  all $i$.   If  we collect together
non-trivial terms with consecutive  indices we obtain a  decomposition
of  this element  as  a product of  signed blocks.   More precisely, a
positive respectively  negative  block  of length   $m+1  \geq 0$   in
$\tilde{W}_n$ is an element of the form
\[s_{i+1}s_{i+2}\cdots s_{i+m} \quad \mbox{ respectively } \quad
t_i's_{i+1} \cdots s_{i+m}\] where  $m,i\geq 0$. If  we denote such an
element     by        $b(i,m)$       then           \[d_1       \cdots
d_n=b(i_1,m_1)b(i_2,m_2)b(i_3,m_3)      \cdots,  \quad  \mbox{ where }
i_2>i_1+m_1, i_3>i_2+m_2 \mbox{ etc}.\] By \cite{GP}, \S 2, each conjugacy 
class in the Coxeter group $W_n$ contains an element in signed block 
form, and such an element is of miminal length in its class if and only 
if all negative blocks are in the beginning, ordered by increasing length. 
In order to reduce an arbitrary element in $W_n$ to such a minimal form, 
it is necessary to interchange by conjugation two consecutive blocks in 
the signed block form of an element in $W_n$. Our aim here is to show 
that similar relations also hold in the braid group.

\addtocounter{thm}{1}
\begin{lem} {\rm (cf.\ \cite{GP}, Proposition~2.4.)}\\
Let  $y=(s_{i+m+1}  \cdots s_{i+1})  (s_{i+m+2} \cdots s_{i+2}) \cdots
(s_{i+m+k+1} \cdots s_{i+k+1})$ for $i,k,m \geq 0$. Geometrically
\ $y$ \ is  a half-twist of  $m+1$ consecutive strands around the next
$k+1$ consecutive strands in the classical braid group.

(a) Let
\begin{eqnarray*} 
w & = & (s_{i+1} s_{i+2} \cdots s_{i+m}) (s_{i+m+2} 
s_{i+m+3} \cdots s_{i+m+k+1}) \quad \mbox{and}\\ 
v & = & (s_{i+1} s_{i+2} \cdots s_{i+k}) (s_{i+k+2} s_{i+k+3} \cdots 
s_{i+k+m+1}).  \end{eqnarray*}
Then $y^{-1}wy=v$  in the braid group.

(b) Let 
\begin{eqnarray*}
w & = & (s_{i+1} s_{i+2} \cdots s_{i+m}) (t_{i+m+1}'s_{i+m+2} 
s_{i+m+3} \cdots s_{i+m+k+1}) \quad \mbox{and} \\
v & = & (t_i's_{i+1} s_{i+2} \cdots s_{i+k}) (s_{i+k+2} s_{i+k+3} \cdots 
s_{i+k+m+1}).  \end{eqnarray*}
Then $y^{-1}wy=v$ in the braid group.

(c) Let \[w = (t_i's_{i+1} s_{i+2} \cdots s_{i+m}) (t_{i+m+1}'s_{i+m+2} 
s_{i+m+3} \cdots s_{i+m+k+1})\] for some $m>k$. Then $y^{-1}wy=v$ in the 
braid group where
\[ v  = (t_i's_{i+1} s_{i+2} \cdots s_{i+k})(t_{i+k+1,i}'s_{i+k+2} s_{i+k+3}
 \cdots s_{i+k+m+1})\] for some $t_{i+k+1,i}'$ (see the proof below).

In (b), (c), the length of $v$ is strictly shorter than the length of $w$.
\end{lem}

\begin{pf} Geometrically, (a), (b) and (c) follow immediately by looking 
at  the corresponding braid  pictures and comparing their closures; note that
in each case we obtain links  of two components.  For an algebraic proof,
we  use the similar relations   in \cite{GP}, Proposition~2.4, on  the
level of $W_n$.   We have to slightly  modify those arguments in order
to derive relations in $\tilde{W}_n$.

(a)  In \cite{GP}, Proposition~2.4(a), it   is shown that the equation
$wy=vy$ holds  in $W_n$ and  that the  expressions  on both sides  are
reduced.  Hence the  left hand side  can be  transformed to  the right
hand side by  a finite sequence of braid  relations.   It follows that
the equation $wy=vy$ also holds in the braid group.

(b) We write $y=(s_{i+m+1} \cdots s_{i+1})y'$. Then 
\[ y^{-1}t_{i+m+1}'y= y^{-1}s_{i+m+1} \cdots s_{i+1}t_i's_{i+1}^{-1} \cdots
s_{i+m+1}^{-1}y = {y'}^{-1}t_i'y'=t_i'.\]  For the last equality, note
that $y'$ is  a product of generators $s_j$  with $j>i+1$ and that all
of these commute with $t_i'$.

Now $(s_{i+1} \cdots s_{i+m})$  commutes with $t_{i+m+1}'$. So we  can
write $w=t_{i+m+1}'w_a$ and $v=t_i'v_a$  where $w_a$ and $v_a$  are as
in part (a) of the Lemma. Hence we deduce, using (a), that
\[ y^{-1}wy=(y^{-1}t_{i+m+1}'y)(y^{-1}w_ay)=t_i'v_a=v.\]

(c) We rearrange the given reduced expression for $y$ as 
\[ y=(s_{i+m+1} \cdots s_{i+m+k+1})  (s_{i+m} \cdots s_{i+m+k}) \cdots  
(s_{i+1} \cdots    s_{i+k+1})\]   and write   $y=y'(s_{i+1}     \cdots
s_{i+k+1})$.  Now note that $y'$ is a product of generators $s_j$ with
$j>i+1$ and  that all of these  commute  with $t_i'$. Then  we compute
that
\begin{eqnarray*}
y^{-1}t_i'y     &          =      &        (s_{i+k+1}^{-1}      \cdots
s_{i+1}^{-1}){y'}^{-1}t_i'y'(s_{i+1}  \cdots    s_{i+k+1}) \\  &   = &
s_{i+k+1}^{-1} \cdots s_{i+1}^{-1}t_i's_{i+1} \cdots s_{i+k+1} \\ & =:
& t_{i+k+1,i}''. \end{eqnarray*} We can  write $w=t_i'w_b$ where $w_b$
is as in part (b). Thus, we deduce, using (b), that
\[y^{-1}wy = (y^{-1}t_i'y)(y^{-1}w_by) = 
(t_{i+k+1,i}''t_i's_{i+1} \cdots s_{i+k})(s_{i+k+2} \cdots s_{i+k+m+1}).\]
Using Lemma~2.3(b), we compute that $t_{i+k+1,i}''t_i'=t_i't_{i+k+1,i+1}''$
where the inverses in $t_{i+k+1,i+1}''$ at all indices bigger than $i+1$ 
are on the left hand side of $t$. A final calculation then shows that 
\[t_{i+k+1,i+1}''(s_{i+1} \cdots s_{i+k})=(s_{i+1} \cdots s_{i+k})
(s_{i+k+1}s_{i+k}^{-1} \cdots s_{i+1}^{-1}t_i's_{i+1} \cdots 
s_{i+k}s_{i+k+1}^{-1}).\]
If we denote the second factor on the right hand by $t_{i+k+1,i}'$ then
$y^{-1}wy=v$ as required.  The proof is complete.  \end{pf}

\begin{rems} \rm 
(a) One can also  show that, if   the conjugation in  (b), (c)  of the
above Lemma is performed step by step (one generator of $y$ at a time)
then this sequence of conjugations can be arranged  in such a way that
the length  of the elements  does  not increase at  each step   (cf.\
\cite{GP}, Proposition~2.4(b,c)).

(b) Consider the elements
\[ w=(s_1) (s_2^{-1}s_1ts_1^{-1}s_2) \quad \mbox{ and } \quad v=ts_2 
\quad \mbox{(in $\tilde{W}_3$)}.\]
If  all the inverses in  $w$ were on the  right  hand side of $t$ then
this element would be  conjugate to $v$  in $\tilde{W}_3$, by part (b)
of the above Lemma  (with $m=1$, $k=i=0$). However,  one can show that
$w$  and $v$ are  not conjugate in $\tilde{W}_n$, and  not even in the
associated Iwahori-Hecke algebra.   This  example indicates  that  the
statements in the above Lemma  are as strong as  possible, and that we
cannot distribute the inverses in some arbitrary order around $t$ when
we want to conjugate elements in block form.

Note, however, that the oriented links obtained by closing the braids
corresponding to the above elements are isotopic (see (5.4) below). In 
particular, the knot invariants constructed in Definition~5.3 must
have the same value on them.
\end{rems}

\section{Trace  functions on the Iwahori-Hecke algebra of type $B$}

In this section, we introduce the  Iwahori-Hecke algebra $H_n$ of type
$B_n$ as a quotient of  the braid group algebra  of $\tilde{W}_n$.  We
also   show how the  main  results  of \cite{GP}  on determining trace
functions on  $H_n$ can be   adapted to  our  present situation  where
reduced expressions  for representatives   of minimal length    in the
conjugacy classes involve some inverses.

\medskip
3.1. Let $A$ be a commutative ring with~$1$ and $Q,q \in A$ two
fixed invertible elements.  Then $H_n$  is an associative algebra over $A$. 
It can be described as a quotient of the
group algebra of the braid  group $\tilde{W}_n$ (over $A$) obtained by
factoring by the ideal generated  by all elements  of the form $t^2 
-(Q-1)t-Q \cdot 1$, $s_i^2-(q-1)s_i- q \cdot 1$ for $i=1,\ldots,n-1$.  We
denote the    image  of $t$  under the   canonical   map $A\tilde{W}_n
\rightarrow H_n$ again by  $t$, and the image of  $s_i$ by  $g_i$, for
all $i$. Then the generators  $t,g_1,\ldots,g_{n-1}$ of $H_n$  satisfy
braid  relations completely analogous to  the braid  relations for the
generators  $t,s_1,\ldots,s_{n-1}$ of   $\tilde{W}_n$.  In
addition, we have the following quadratic relations.
\[ t^2=(Q-1)t+Q \cdot 1  \quad \mbox{ and } \quad g_i^2=(q-1)g_i+q \cdot 1 
\mbox{ for all $i$}.\]
Let $w \in W_n$ and assume that we  are given a reduced expression for
$w$ as a  product    of generators $t,s_1,\ldots,s_{n-1}$.  Then   the
corresponding element      of $H_n$  in  terms     of  the  generators
$t,g_1,\ldots,g_{n-1}$ is independent of the chosen reduced expression
for $w$.  We may therefore  denote this element in $H_n$ unambiguously
by   $g_w$. It is  known that  the set  of elements $\{g_w  \mid w \in
W_n\}$  forms   an $A$-basis of   $H_n$.  (For   all these  facts see
\cite{Bour},  Chap.\  IV,  \S 2, Ex.\   23.)   We then also   have the
following relations.
\[ g_w g_{w'}=g_{ww'} \quad \mbox{ if $l(ww')=l(w)+l(w')$}.\]
The embeddings $W_0 \subset W_1 \subset W_2 \subset \ldots$ of (2.1)
induce corresponding embeddings of algebras $H_0 \subset H_1 \subset H_2 
\subset \ldots$ and we shall denote
\[ H:=\bigcup_{n \geq 0} H_n.\]

\medskip
3.2.  The fact that $q$ is invertible in $A$ implies that the generators
$g_i$ are also invertible in $H$. In fact, we have that
\[ g_i^{-1}=q^{-1}g_i+(q^{-1}-1) \cdot 1 \in H_n.\]
Thus, the images of the elements $t_i',t_{i,j}' \in \tilde{W}_n$ under
the map  $A\tilde{W}_n  \rightarrow H_n$ are  well-defined elements in
$H_n$, and we shall denote them by the same symbols. We also write
\[{\cal D}_n' = \{ 1,g_{n-1}, t_{n-1}' \}\]
and,   similarly, for  ${\cal   R}_n'$, ${\cal  D}_{n,i}'$ and  ${\cal
R}_{n,i}'$ (cf.\ (2.2)). With   these conventions, all results   about
commutation and conjugation of the various special elements considered
in the previous section carry over without change to $H_n$.

Let $w \in W_n$ and write $w=r_1 \cdots r_n$ with $r_i \in {\cal R}_i$
for  all $i$. Since this expression is reduced we also have $g_w=
g_{r_1} \cdots g_{r_n}$. For  each $r_i$ let   $r_i' \in {\cal  R}_i'$  be 
the corresponding element in $H_n$ (where the $s_j$ are replaced by $g_j$,
and $t_j$  by   $t_j'$).  Let $n_w  \geq  0$  be the   total number of
inverses in the   terms $r_1',\ldots,r_n'$. Using  the above inversion
formula it then follows that
\[ g_w=q^{n_w}r_1' \cdots r_n'+ \mbox{$A$-linear combination of elements
$g_v$ with $l(v) < l(w)$}.\] One consequence of this  is the fact that
the elements $\{r_1' \cdots r_n' \mid r_i' \in  {\cal R}_i'\}$ form an
$A$-basis of   $H_n$,  and if   we order  the   elements of  $W_n$  by
increasing  length then the matrix  performing the base  change to the
old  basis $\{ g_{r_1    \cdots r_n} \mid   r_i  \in {\cal  R}_i\}$ is
triangular with powers of $q$ along the diagonal.

\medskip 
3.3. Let $\{C\}$ be the set of conjugacy classes of $W_n$ and let $w_C$ be 
an element of minimal length in $C$ which admits a decomposition as a 
product of negative blocks (ordered by increasing length) followed
by various positive blocks (see (2.4)). We can even fix a unique choice
of $w_C$ if we also require that the positive blocks  have increasing 
length. We then define an element $g_C \in  H_n$ by taking an
expression $w_C=d_1 \cdots d_n$ (with $d_i \in {\cal D}_i$) and replacing 
each  $s_i$ by  $g_i$ and each $t_i$ by $t_i'$.  As in (3.2) we have (with 
$n_C:=n_{w_C}$) that
\[ g_{w_C} =q^{n_C}g_{C}+   \mbox{$A$-linear combination of elements $g_w$ 
with $l(w)<l(w_C)$}.\] 
A  trace  function   on $H_n$   is  an  $A$-linear  map   $\varphi:H_n
\rightarrow  A$ such that  $\varphi(h h')=\varphi(h' h)$ for all $h,h'
\in   H_n$.  By \cite{GP}, each trace   function  on $H_n$ is uniquely
determined  by  its   values   on  the  elements   $g_{w_C}$,  for all
$C$. Conversely, given  a set  of elements  $a_C \in A$,  one for each
conjugacy class $C$, there exists a unique trace function $\varphi$ on
$H_n$ such that $\varphi(g_{w_C})=a_C$ for all~$C$.  Using the above 
relations, we deduce that these results on the determination of trace 
functions  remain valid when we replace each $g_{w_C}$ by $g_C$, for all~$C$.

The following result will show how to reduce the computation of the value
of a  trace function  on any element to the values on elements in signed
block form.

\addtocounter{thm}{3}
\begin{prop}
For each $h \in  H_n$ there exists a  finite (non-empty) subset  $I(h)
\subseteq   A  \times  {\cal  D}_{1,0}'  \times   \ldots \times  {\cal
D}_{n,0}'$ such that \[\varphi(h) =   \sum_{(r,d_1^{*},\ldots,d_n^{*})
\in I(h)} r \varphi(d_1^{*} \cdots d_n^{*}),\] for all trace functions
$\varphi$ on $H_n$.
\end{prop}

\begin{pf}
The result clearly holds  if $n=1$. Now  let  $1<j \leq n$ and  assume
that we have already found  a finite (non-empty) subset $I_j \subseteq
H_j \times D_{j+1,j-1}' \times  \ldots \times {\cal D}_{n,j-1}'$  such
that \[ \varphi(h)= \sum_{(h_j,d_{j+1}',\ldots,d_n') \in I_j} 
\varphi(h_jd_{j+1}' \cdots d_{n}'),\]
for  all trace  functions $\varphi$  on    $H_n$. We will  proceed  by
downward induction on $j$. For $j=n$ there is nothing to prove. We now
show  how to   obtain an  analogous   statement with  $j$ replaced  by
$j-1$. This is done as follows.

Consider one element  $(h_j,d_{j+1}',  \ldots,d_n') \in  I_j$. The
element $h_j$ is an $A$-linear combination of basis elements $g_w$
with $w \in W_j$. By (3.2), this can be rewritten as an $A$-linear
combination of products $r_1' \cdots r_j'$ with $r_i' \in {\cal R}_i'$.
Collecting terms with a fixed value of $r_j'$, we obtain a finite 
(non-empty) subset   $R(h_j) \subseteq H_{j-1}  \times {\cal  R}_j'$   
such  that \[h_j=\sum_{(h_{j-1},r_j')  \in R(h_j)} h_{j-1}r_j'.\] 
We  can write $r_j'=d_j'r_{j-1}''$   with  $d_j' \in  {\cal D}_j'$  
and $r_{j-1}'' \in {\cal R}_{j-1}'$ (cf.\  (2.1)). Now  the element 
$r_{j-1}''$ either is a product of various generators $g_1,\ldots,g_{j-2}$
or is like the element considered in Lemma~2.3(c). In any case,  it 
commutes with $d_{j+1}'$ up to (possibly) changing some inverses at 
indices bigger than $j-2$, and similarly for $d_{j+2}',\ldots,d_{n}'$. 
Thus, we conclude that 
\[ h_jd_{j+1}'\cdots d_n' = \sum_{(h_{j-1},r_j') \in R(h_j)} 
h_{j-1}d_j'd_{j+1}'' \cdots d_n''r_{j-1}'',\] where $d_{j+1}'' \in   
{\cal D}_{j+1,j-2}',\ldots,d_n'' \in {\cal D}_{n,j-2}'$. 
So we have that \[\varphi(h)=\sum_{(h_j,d_{j+1}',\ldots,d_n') \in I_j}
\sum_{(h_{j-1},r_j') \in R(h_j)} \varphi(r_{j-1}''h_{j-1}d_j'd_{j+1}''
\cdots d_n'').\]  We  can combine the index   sets to a  new  set  
$I_{j-1}  \subseteq H_{j-1}  \times {\cal D}_{j,j-2}' \times \ldots 
\times  {\cal D}_{n,j-2}'$ and arrive  at a new  expression  as above  
with  $j$ replaced  $j-1$.  We  repeat this process until  we  arrive 
at  $j=1$.  Then $H_1=\langle 1,  t \rangle= \langle {\cal D}_1' 
\rangle$, and we are done.  \end{pf}

In fact, the above proof also yields the following extension.

\begin{cor}
Let $h \in H_n$ and $d_{n+1} \in {\cal D}_{n+1}', \ldots, d_{n+m} \in 
{\cal D}_{n+m}'$ (for some $m \geq 0$). Then, for all trace functions
$\varphi$ on $H_n$, we have that
\[\varphi(hd_{n+1} \cdots d_{n+m})=\sum_{(r,d_1^{*},\ldots,d_n^{*}) \in I(h)} 
r\varphi(d_1^{*} \cdots d_n^{*}d_{n+1}^{*} \cdots d_{n+m}^{*}),\]
for some elements $d_{n+1}^{*} \in {\cal D}_{n+1,0}',\ldots,
d_{n+m}^{*} \in {\cal D}_{n+m,0}'$ (depending on the various elements in 
$I(h)$).
\end{cor}

\section{Markov traces for Iwahori-Hecke algebras of type $B$}

Jones writes in \cite{Jones}, p.346, that there should be analogues of
Ocneanu's trace for  Iwahori-Hecke  algebras other than those  of type
$A$. The  trace given  in \cite{L2} was  the  first such  analogue for
$B$-type  Iwahori-Hecke   algebras.  The  aim of  this   section is to
classify  {\it all} such  `Markov' traces on Iwahori-Hecke algebras of
type $B$, based on the results in the previous sections.

\begin{defn} 
Let $z \in A$ and $\tau:H \rightarrow A$ be an $A$-linear map. Then $\tau$
is called a Markov trace (with parameter $z$) if the following conditions 
are satisfied.
\begin{itemize}
\item[(1)] $\tau$ is a trace function on $H$.
\item[(2)] $\tau(1)=1$ (normalization).
\item[(3)] $\tau(hg_n)=z\tau(h)$ for all $n \geq 1$ and $h \in H_n$.
\end{itemize}
\end{defn}

\noindent We note that all generators $g_i$ (for $i=1,2,\ldots$) are 
conjugate in $H$. In particular, any trace function on $H$ must have the 
same value on these elements. This explains why the parameter $z$ is 
independent of $n$ in rule (3) of this definition. 

Let   us  consider    the  subalgebra   $H'$ of    $H$  generated   by
$g_1,g_2,\ldots$.   Then  $H'$ is  the algebra  considered by Jones in
\cite{Jones}, \S 5 (infinite  union over all Iwahori-Hecke algebras of
type  $A_n$ with parameter  $q$).   Moreover, the restriction of  any
Markov trace $\tau$ on $H$ to $H'$  yields Ocneanu's original trace as
in  \cite{Jones},  Theorem~5.1,  which is uniquely   determined by the
parameter~$z$.

The next result describes a set of elements in $H$ which is sufficient
to determine a  Markov trace $\tau$. We shall see that this set
is in fact as small as possible.

\begin{lem}
Let $\tau:H \rightarrow A$ be a Markov trace (with parameter~$z \in A$).
\begin{itemize}
\item[(a)] If $n \geq 1$, $m \geq 0$ and $h \in H_n$, then 
\begin{eqnarray*}
 \tau(hg_nt_{n+1}' \cdots t_{n+m}') & = & z \tau(ht_n' \cdots 
t_{n+m-1}') \quad \mbox{ and }\\
 \tau(ht_{n+1}' \cdots t_{n+m}') & = &  \tau(ht_n' \cdots 
t_{n+m-1}') \end{eqnarray*}
\item[(b)] If $h=d_1 \cdots d_n \in H$ where $d_i \in {\cal D}_i'$ for all
$i$ then  \[ \tau(h)=z^{a(h)} \tau(t_0't_1' \cdots  t_{b(h)-1}')\] where
$a(h)$  is    the   number  of    $d_i$    which   are  in  the    set
$\{g_1,\ldots,g_{n-1}\}$  and $b(h)$ is the number  of $d_i$ which are
conjugate to $t$.
\item[(c)] $\tau$ is uniquely determined by its values on the elements 
in the set \[ \{ t_0' t_1' \cdots t_{k-1}' \mid k=1,2,\ldots \}.\]
\end{itemize}
\end{lem}

\begin{pf} To prove the first relation in (a) we will proceed by induction 
on $m$.    If $m=0$   then we    can   apply directly    rule (3)   in
Definition~4.1.  Now let  us assume that  $m>0$.  We  have to evaluate
the expression
\[\tau(hg_nt_{n+1}' \cdots t_{n+m}').\]
We write $t_{n+1}'=g_{n+1}t_n'g_{n+1}^{-1}$ and observe that $g_{n+1}^{-1}$
commutes with $t_{n+2}',\ldots,t_{n+m}'$ by Lemma~2.3(a). Since $\tau$ is
a trace our expression is equal to
\[\tau(g_{n+1}^{-1}hg_ng_{n+1}t_n't_{n+2}' \cdots t_{n+m}').\]
Now $h$   lies in $H_n$,  that is,  $h$ only   involves the generators
$t,g_1,\ldots,g_{n-1}$.     It follows    that    $h$   commutes  with
$g_{n+1}^{-1}$.    Using      moreover        the   braid     relation
$g_{n+1}^{-1}g_ng_{n+1}=g_ng_{n+1}g_n^{-1}$, the above  expression can
be rewritten as
\[\tau(hg_ng_{n+1}g_n^{-1}t_n't_{n+2}' \cdots t_{n+m}').\] 
If we write $t_n'=g_nt_{n-1}'g_n^{-1}$,  the left hand term $g_n$ will
cancel and then $g_{n+1}$ commutes with $t_{n-1}'$. Now our expression
reads
\[\tau(hg_nt_{n-1}'g_{n+1}g_n^{-1}t_{n+2}' \cdots t_{n+m}').\] 
The element $g_n^{-1}$ commutes with all terms to the right of it. Hence our
expression is equal to
\[\tau(g_n^{-1}hg_nt_{n-1}'g_{n+1}t_{n+2}' \cdots t_{n+m}').\]
We write $h':=g_n^{-1}hg_nt_{n-1}'$ and observe that this element lies
in $H_{n+1}$. So we can apply the induction and obtain that
\[\tau(hg_nt_{n+1}' \cdots t_{n+m}')=\tau(h'g_{n+1}t_{n+2}' \cdots t_{n+m}')=
z \tau(h't_{n+1}' \cdots t_{n+m-1}').\]
We insert the expression for $h'$ again, note that $g_n^{-1}$ commutes
with $t_{n+1}' \cdots t_{n+m-1}'$, and conclude that
\[\tau(h't_{n+1}' \cdots t_{n+m}')= 
\tau(hg_nt_{n-1}'g_n^{-1} t_{n+1}' \cdots t_{n+m-1}')=
\tau(ht_n't_{n+1}' \cdots t_{n+m-1}').\]
Putting  things together we see that  this completes  the proof of the
first relation. The proof  of the second is  achieved by  an analogous
computation  (with $g_n$ replaced by~$1$).  In  order  to prove (b) we
consider an element
\[ h=d_1 \cdots d_n \in H_n \quad \mbox{ where $d_i \in {\cal D}_i'$ for all
$i$.} \] Using (a) and induction on $n$ we deduce that 
\[ \tau(h)=z^{a(h)} \tau(t_0't_1' \cdots t_{b(h)-1}')\]
where $a(h),b(h)$ are defined as above. Using (3.3) we 
conclude that these equations determine $\tau$ uniquely, proving (c).  
\end{pf}

\begin{thm}
Let $z,y_1,y_2 \ldots \in A$. Then there exists a unique Markov trace 
$\tau$ on $H$ with parameter $z$ such that 
\[ \tau(t_0't_1't_2' \cdots t_{k-1}')=y_k \quad \mbox{ for all $k \geq 1$.}\]
\end{thm}

\begin{pf} Uniqueness was already proved in Lemma~4.2. Using the facts 
summarized  in  (3.3) we can  define  a  trace function  $\tau$ on $H$
satisfying $\tau(1)=1$ and \[ \tau(g_C)=z^ay_b\] where $g_C=d_1 \cdots
d_n$   with   $d_i  \in  {\cal D}_i'$   and   the elements $a=a(g_C)$,
$b=b(g_C)$ are defined as in Lemma~4.2. Thus, the existence of a trace
function  $\tau$    on   $H$   satisfying   conditions~(1),  (2)    in
Definition~4.1 is already established. The problem is to show that (3)
holds.

This will  be done by an  induction, as follows.   For $N \geq  0$ let
$H(\leq N)$ be the $A$-subspace of $H$ generated by all elements $g_w$
with $w \in W$ and $l(w) \leq N$. We  shall prove the following claim,
for all $N \geq 0$.
\begin{itemize}
\item[(*)] {\em Let $h \in H_n$ and $d_{n+i}^{*} \in  
{\cal D}_{n+i,0}'$ (for some $n,m \geq 1$ and $i=1,\ldots,m$) such that 
$hd_{n+1}^{*}\cdots d_{n+m}^{*} \in H(\leq N)$.
Then} 
\begin{eqnarray*}
\tau(hd_{n+1}^{*} \cdots d_{n+m}^{*}) & = & \tau(hd_{n+1} \cdots 
d_{n+m}) \quad \mbox{ and} \\
\tau(hg_nd_{n+2}^{*}d_{n+3}^{*} \cdots d_{n+m}^{*}) & = & z\tau(hd_{n+2}^{*} 
\cdots d_{n+m}^{*}). \end{eqnarray*} 
\end{itemize}
(Here, we used  the following convention. For each  $i$, we  denote by
$d_i$ the analogous element as $d_i^{*}$ with the inverses (if any) on
the right hand side of $t$.) If this is proved for all $N$ then $\tau$
will   satisfy condition  (3)    in    Definition~4.1 as a     special
case. Clearly, (*) is true  for $N=0$. Now  let $N>0$. We proceed in a
number of steps.

\medskip
{\em Step~1.} At first we show that $\tau(hd_{n+1}^{*} \cdots d_{n+m}^{*})
=\tau(hd_{n+1} \cdots d_{n+m})$. By induction on $m$ we may assume
that $d_{n+2}^{*}=d_{n+2},\ldots, d_{n+m}^{*}=d_{n+m}$. We can also assume
that
\[d_{n+1}^{*}=T_{n+1,i}:=g_{n} \cdots g_{i+1}g_i^{\mp 1}t_{i-1,0}'
g_i^{\pm 1}g_{i+1}^{-1} \cdots g_{n}^{-1} \quad \mbox{ for some  $i \leq n$}.
\]
(That is, the inverses are already fixed at  indices   bigger than $i$.)    
Then $T_{n+1,i-1}$ is  the analogous element where  the inverse at index  
$i$ is also  correct.  Now assume that the sign in $g_i^{\mp 1}$ is $-1$.
We will show that $\tau(hT_{n+1,i}d_{n+2} \cdots d_{n+m})=\tau(hT_{n+1,i-1}
d_{n+2} \cdots d_{n+m})$.  This is done as follows.  In $T_{n+1,i}$  and 
$T_{n+1,i-1}$ we  replace $g_i^{-1}$ by $(q^{-1}-1)g_i+q^{-1} \cdot 1$. 
Then $hT_{n+1,i}d_{n+2} \cdots d_{n+m}=(q^{-1}-1)S_1+ q^{-1}S_2$ and 
$hT_{n+1,i-1}d_{n+2} \cdots d_{n+m}=(q^{-1}-1)S_1+q^{-1}S_3$ where
\begin{eqnarray*}
S_1 & = & h (g_n \cdots g_{i+1}g_i t_{i-1,0}'g_ig_{i+1}^{-1} \cdots g_n^{-1})
                    d_{n+2} \cdots d_{n+m},\\
S_2 & = & h (g_n \cdots g_{i+1}t_{i-1,0}'g_ig_{i+1}^{-1} \cdots g_n^{-1}) 
                    d_{n+2} \cdots d_{n+m},\\
S_3 & = & h (g_n \cdots g_{i+1}g_it_{i-1,0}'g_{i+1}^{-1} \cdots g_n^{-1})
                    d_{n+2} \cdots d_{n+m}.
\end{eqnarray*}
Note that $S_2,S_3 \in H(\leq N-1)$ so that we can use our inductive
hypotheses in the evaluation of $\tau$ on these elements.
Let us first consider $S_2$. At first, we note that $g_n \cdots g_{i+1}$
commutes with $t_{i-1,0}'$. So we obtain that
\begin{eqnarray*}
\tau(S_2) & = & \tau(ht_{i-1,0}'(g_n \cdots g_{i+1}g_ig_{i+1}^{-1} 
               \cdots g_n^{-1})d_{n+2} \cdots d_{n+m}) \\
    & = & \tau(ht_{i-1,0}' (g_i^{-1} \cdots g_{n-1}^{-1}g_ng_{n-1} 
                    \cdots g_i)d_{n+2} \cdots d_{n+m}).\end{eqnarray*}
Now $g_{n-1} \cdots g_{i}$ commutes with all terms to the right of it,
by Lemma~2.3(a). Hence we conclude that
\[ \tau(S_2)= \tau((g_{n-1} \cdots g_i)ht_{i-1,0}' (g_i^{-1} 
         \cdots g_{n-1}^{-1})g_nd_{n+2} \cdots d_{n+m}). \]
We write the argument of $\tau$ as $h'g_nd_{n+2} \cdots 
d_{n+m}$ with $h' \in H_n$. Since this element lies in $H(\leq N-1)$
we can use induction to deduce that 
\[\tau(S_2)=\tau(h'g_nd_{n+2} \cdots d_{n+m})= 
z\tau(h'd_{n+2} \cdots d_{n+m}).\]
Now we observe that $h'$ is a conjugate of $ht_{i-1,0}'$ and the conjugating
element $g_{n-1} \cdots g_i$ commutes with $d_{n+2} \cdots d_{n+m}$, again
using Lemma~2.3(a). Thus, we finally compute that
\[\tau(S_2) = z \tau(ht_{i-1,0}'d_{n+2} \cdots d_{n+m}).\]
Now we follow the same procedure with $S_3$ and find that
\begin{eqnarray*}
\tau(S_3) & = & \tau(h (g_n \cdots g_{i+1}g_ig_{i+1}^{-1} \cdots g_n^{-1})
                                    t_{i-1,0}'d_{n+2} \cdots d_{n+m}) \\
    & = & \tau(h (g_i^{-1}  \cdots g_{n-1}^{-1}g_ng_{n-1} 
                            \cdots g_i)t_{i-1,0}'d_{n+2} \cdots d_{n+m}).
\end{eqnarray*}
Now we  would like  to  interchange $t_{i-1,0}'$  and $d_{n+2}  \cdots
d_{n+m}$ but  this   is   not  necessarily possible   on   the algebra
level.  However, it still works  modulo the kernel  of  $\tau$, as the
following   argument shows.   We  write  $t_{i-1,0}'=x^{-1}tx$   where
$x=g_1^{\pm 1} \cdots g_{i-1}^{\pm  1}$.  By Lemma~2.3(a), the element
$x$ will commute with $d_{n+2} \cdots d_{n+m}$.  The next factor, $t$,
only commutes with $d_{n+2}  \cdots d_{n+m}$ up to (possibly) changing
some  inverses. Thus, the  term  on the right   hand side of the above
equality is equal to
\[ \tau(tx h (g_i^{-1} \cdots g_{n-1}^{-1}g_ng_{n-1} 
                         \cdots g_i)x^{-1}d_{n+2}' \cdots d_{n+m}'),\]
for some $d_{n+2}' \in {\cal D}_{n+2,0}',\ldots, d_{n+m}' 
\in {\cal D}_{n+m,0}'$. Now we can apply our inductive hypothesis to 
conclude that this equals
\[ \tau(tx h (g_i^{-1} \cdots g_{n-1}^{-1}g_ng_{n-1} 
                         \cdots g_i)x^{-1}d_{n+2} \cdots d_{n+m}).\]
Now, again by Lemma~2.3(a), the element $x^{-1}$ will commute with the
terms to the right of it. So, finally, we obtain that
\[ \tau(S_3) = \tau(t_{i-1,0}'h (g_i^{-1} \cdots g_{n-1}^{-1}g_ng_{n-1} 
                         \cdots g_i)d_{n+2} \cdots d_{n+m}).\]
Once more, we use Lemma~2.3(a) to conclude that $g_{n-1} \cdots g_i$ 
commutes with all terms to the right of it. So we arrive at the equality
\[ \tau(S_3)=\tau((g_{n-1} \cdots g_i)t_{i-1,0}'h 
(g_i^{-1} \cdots g_{n-1}^{-1}) g_nd_{n+2} \cdots d_{n+m}).\]
In this situation, we can argue similarly as in the evaluation
of $\tau(S_2)$. This evaluation will result in  an analogous expression as 
before, but with the terms $h$ and $t_{i-1,0}'$ interchanged. Thus, we
conclude that
\[\tau(S_3)=z\tau(t_{i-1,0}'hd_{n+2} \cdots d_{n+m}).\] Arguing as above,
we see that $t_{i-1,0}'$ and $d_{n+2} \cdots d_{n+m}$ can be 
interchanged modulo the kernel of $\tau$.  So, eventually, we find that
$\tau(S_2)=\tau(S_3)$ and, hence, that $\tau(hT_{n+1,i}d_{n+2} \cdots
d_{n+m})=\tau(hT_{n+1,i-1}d_{n+2} \cdots d_{n+m})$, as required. The proof 
of the assertion of Step~1 is now completed by repeating the whole process 
with $T_{n+1,i-1}$ and so on, until all inverses are fixed.

\medskip
{\em  Step~2.} Now  we  consider  the   case where $h=d_1^{*}   \ldots
d_n^{*}$, for some $d_i^{*} \in {\cal D}_{i,0}'$. As before, let $d_i$
be the  element of ${\cal  D}_i'$ corresponding to $d_i^{*}$ where all
inverses  (if any) are  in the right position.  Using Step~1, we first
see   that   $\tau(d_1^{*}    \cdots  d_{n+m}^{*})=   \tau(d_1  \cdots
d_{n+m})$. Now, we will prove that
\[\tau(d_1^{*}\cdots d_{n+m}^{*})=z^a y_{b},\] where 
$a=a(d_1 \cdots d_{n+m})$  and $b=b(d_1 \cdots d_{n+m})$. (Note  that,
in the  definition  of $a,b$ in Lemma~4.2(b),   it does not   matter  
whether we  take $d_i^{*}$ or $d_i$.)  As explained in (2.4),  (3.3), the  
element $d_1 \cdots d_{n+m}$ can be regarded as a  product of positive and 
negative blocks. If this element is equal to $g_C$ for some conjugacy class 
$C$ then we are done by the  defining equation for $\tau$.  If this is not
the case, then some positive block is followed by a negative block or
some  negative    block is followed  by   a  strictly shorter negative
block. We   can  then  use Lemma~2.5 to     conjugate our  element  to
${d_1'} \cdots {d_{n+m}'}  \in H_n(\leq N-1)$ with ${d_i'}
\in  {\cal D}_{i,0}'$ and  with   the same  signed block structure  as
before, that is,  the values of  $a$ and $b$ haven't changed after this
conjugation. Thus, we conclude that
\[ \tau(d_1^{*} \cdots d_{n+m}^{*}) = \tau(d_1 \cdots d_{n+m}) 
  = \tau({d_1'} \cdots  {d_{n+m}'}) = z^ay_b,\] where the last
equality is by induction on $N$.  Our claim is proved.

\medskip
{\em Step 3.} Finally, let $h \in H_n$ and
$d_{n+i}^{*} \in {\cal D}_{n+i,0}'$ for $2 \leq i \leq m$. Now, we show
that \[ \tau(hg_nd_{n+2}^{*}d_{n+3}^{*} \cdots d_{n+m}^{*}) = 
z\tau(hd_{n+2}^{*} \cdots d_{n+m}^{*}). \]
Using Step~1, we have that $\tau(hg_nd_{n+2}^{*} \cdots d_{n+m}^{*})=
\tau(hg_nd_{n+2} \cdots d_{n+m})$ where $d_{n+i} \in {\cal D}_{n+i}'$ 
corresponds to $d_{n+i}^{*}$ as above. Using Corollary~3.5, the latter 
trace is equal to
\[\sum_{(r,d_1^{*},\ldots,d_{n}^{*}) \in I(h)} r\tau(d_1^{*}\cdots d_n^{*}
g_nd_{n+2}^{*}  \cdots  d_{n+m}^{*}),\]  for  some $d_{n+2}^{*} \in {\cal
D}_{n+2,0}',\ldots,d_{n+m}^{*} \in  {\cal D}_{n+m,0}'$ (depending   on
the various elements in $I(h)$).   Let  us consider  one term in  this
sum,  corresponding to   $(r,d_1^{*},  \ldots,d_n^{*}) \in   I(h)$. 
We shall write $d_1^{*} \cdots   d_n^{*}  
g_nd_{n+2}^{*}  \cdots  d_{n+m}^{*}$  in   the  form
$h_1g_nh_2$. Let $a=a(h_1g_nh_2)$, $b=b(h_1g_nh_2)$  and
$a'=a(h_1h_2)$, $b'=b(h_1h_2)$.  Then, clearly,  $a'=a-1$  and $b=b'$.
Using this and Step~2, the value of $\tau$ on our element is given by
\[ rz^ay_{b}=rz^{a'}y_{b'}z.\]
Thus, for each term in the above sum, we obtain a factor $z$ as the
expense of cancelling the factor $g_n$ in that term. We conclude that
\[ \tau(hg_nd_{n+2}^{*} \cdots d_{n+m}^{*})=
z\tau(hd_{n+2}^{*} \cdots d_{n+m}^{*}).\]
The proof is complete.  \end{pf}

\medskip
4.4. In the course of the above proof, we have shown the following remarkable
property of a Markov trace $\tau$.

\medskip
\noindent {\em  Let $h \in H_n$  and $d_{n+1}^{*}  \in {\cal D}_{n+1,0}',
\ldots, d_{n+m}^{*} \in {\cal D}_{n+m,0}'$ for some $m \geq 1$. Then
\[ \tau(hd_{n+1}^{*} \cdots d_{n+m}^{*})=\tau(hd_{n+1} \cdots d_{n+m}),\] 
where, as before, $d_i$ is the element in ${\cal D}_i'$ corresponding 
to $d_i^{*}$.}
\medskip

For $h \in H_n$ let $I(h)$ be the corresponding set given by
Proposition~3.4. For each $i=(r,d_1^{*},\ldots,d_n^{*}) \in I(h)$ we
write $r(i)=r$, $a(i)=a(d_1 \cdots d_n)$ and $b(i)=b(d_1 \cdots d_n)$. Then
the above property in combination with Proposition~3.4 and Lemma~4.2 yields 
the following rule for computing $\tau(h)$.
\[ \tau(h)=\sum_{i \in I(h)} r(i) z^{a(i)}y_{b(i)}.\]
Note that an algorithm for computing $I(h)$ is given by the inductive
proof of Proposition~3.4.

\addtocounter{thm}{1}
\begin{prop}
Let $z,y \in A$ and $\tau:H \rightarrow A$ be a Markov trace with
parameter $z$ and such that $\tau(t_0't_1' \cdots t_{k-1}')=y^{k}$ 
for all $k \geq 1$. Then \[ \tau(ht_{n,0}')=y\tau(h) \quad 
\mbox{ for all $n \geq 0$ and $h \in H_n$}.\]
\end{prop}

\begin{pf} Following the steps in the proof of Theorem~4.3 we see that
it is sufficient to consider the case where $h=d_1^{*} \cdots d_n^{*}$
with $d_i^{*} \in {\cal D}_{i,0}'$ for all~$i$. The result in this
case follows from Lemma~4.2. \end{pf}

\begin{rems}{\rm  (a) Given elements $z,y  \in  A$, the  existence of a  Markov 
trace
$\tau= \tau_{z,y}$ as in the previous Proposition can also be proved
along the lines of the approach followed by  Jones \cite{Jones}, \S 5,
in the proof of Ocneanu's  Theorem for Iwahori-Hecke algebras of  type
$A_n$. Such a  proof is sketched in \cite{L2}  (see also \cite{L1}, \S
3.3). It is based on the observation that, for all $n \geq 0$, the map
\[ C_n: H_n \oplus H_n \oplus H_n \otimes_{H_{n-1}} H_n \longrightarrow 
H_{n+1}, \quad  a+b+c \otimes  d  \mapsto a+bt_n'+cg_nd \] defines  an
isomorphism  of $(H_n,H_n)$-bimodules.   Analogously to  the  proof of
\cite{Jones}, (5.1),  we  can now   define  an $A$-linear  map  $\tau$
inductively   on  $H=\bigcup_{n}   H_n$   by  the rules:  $\tau(1)=1$,
$\tau(bt_n')=y\tau(b)$  and $\tau(cg_nd)=z\tau(cd)$  where  $b,c,d \in
H_n$  and $n \geq  0$.  It   follows easily  that  this map  satisfies
conditions (2) and  (3) of Definition 4.1, and  the problem then is to 
show that (1) holds, that is, to show that $\tau$ is a trace function. 
This verification is a lengthy and tedious calculation  that we do not 
want to reproduce here. We do not see, however,  how this method could 
be modified so as to give an alternative  proof  of  our more  general
Theorem~4.3, too. Yet another construction of the trace in \cite{L2} was given by 
T.~tom
Dieck in \cite{tD} using Turaev's $R$--matrix approach (see \cite{T2}). It would be
interesting to find such an $R$--matrix interpretation of our more general traces in
Theorem 4.3, too.

\medskip
\noindent (b) For the definition of Markov traces it doesn't matter
whether we use the elements $t_i$ or $t_i'$ in the signed block form
of elements. From Theorem~4.3 we see that this only plays a role
in the formulation of the {\em initial conditions} determining the trace,
and it would not be clear how to do this in terms of the elements $t_i$.
This gives an explanation for using $t_i'$ rather than $t_i$.}
\end{rems}

4.7. We can use the results of this section to obtain a classification
of Markov traces for Iwahori--Hecke algebras of type~$D$, in the following
way. First note that if we define $u:=ts_1t$ then the elements $u,s_1,
\ldots,s_{n-1}$ generate a subgroup $W_n' \subset W_n$ which is the
finite Coxeter group of type~$D_n$ with relations given by the following diagram.

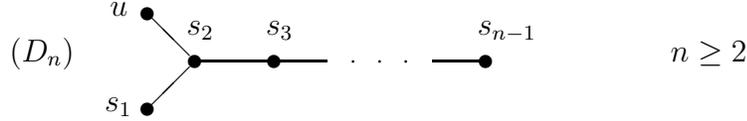
\begin{picture}(300,50)
\put(  0, 20){$(D_{n})$}
\put( 52, 38){\circle*{5}}
\put( 52,  2){\circle*{5}}
\put( 52, 38){\line(1,-1){16}}
\put( 52,  2){\line(1,1){16}}
\put( 70, 20){\circle*{5}}
\put( 70, 20){\line(1,0){30}}
\put( 100, 20){\circle*{5}}
\put( 100, 20){\line(1,0){20}}
\put(130, 20){\circle*{1}}
\put(140, 20){\circle*{1}}
\put(150, 20){\circle*{1}}
\put(160, 20){\line(1,0){20}}
\put(180, 20){\circle*{5}}
\put( 38, 37){$u$}
\put( 36,  1){$s_1$}
\put( 67, 30){$s_2$}
\put( 97, 30){$s_3$}
\put(177, 30){$s_{n-1}$}
\put(250, 20){$n \geq 2$}
\end{picture}

\noindent We shall use the convention that $W_1'=\{1\}$.

The above embedding also works on the level of the Iwahori--Hecke algebras if
we let $Q=1$. Indeed, denoting $u:=tg_1t \in H_n$ in this case, we compute
that $u^2=(q-1)u+q \cdot 1$. Thus, the elements $u,g_1,\ldots,g_{n-1}$
generate a subalgebra $H_n' \subset H_n$ which is the Iwahori--Hecke algebra
of type~$D_n$. As in the $B$-type case, we have embeddings $W_1' \subset W_2'
\subset \ldots$ and $H_1' \subset H_2' \subset \ldots$, and we denote $H':=
\bigcup_{n \geq 1} H_n'$. In analogy to Definition~4.1 we say that a trace 
function $\tau:H' \rightarrow A$ is a Markov trace with parameter $z \in A$ if
$\tau(1)=1$ and $\tau(hg_n)=z\tau(h)$ for all $n \geq 1$ and $h \in H_n'$. Now 
we can state (for the proof see \cite{G}, Section~6):
\begin{itemize}
\item[(a)] {\em Every Markov trace on $H'$ is the restriction of a Markov
trace on~$H$ (with the same parameter~$z \in A$).}
\item[(b)] {\em A Markov trace $\tau$ on $H'$ (with parameter~$z \in A$) is
uniquely determined by its values on the elements in the set 
\[ \{u_1' \cdots u_{2k-1}' \mid k=1,2,\ldots \}\]
where $u_i':=g_i \cdots g_2ug_1^{-1}g_2^{-1} \cdots g_i^{-1}$ for all~$i \geq 1$.}
\end{itemize}
The elements $u_i'$ play an analogous role as the elements $t_i'$ in the 
$B$-type case. Note that under the above embedding $H' \subset H$ we have
$u_i'=tt_i'$ for all $i \geq 1$.

\section{The knot invariants related to the Hecke algebras of $B$-type}

5.1.  Knots  and links  inside  a solid  torus $T$ can be represented  
unambiguously  by  `mixed' knots/links   in  $S^3$  which contain one 
oriented,  unknotted, pointwise fixed  component (the core of the 
complementary unknotted solid torus in $S^3$).  An example of a
mixed link is illustrated below.

$$\vbox{\picill1.7inby1in(mein3)  }$$ 

\noindent So, two links  $L_1$, $L_2$ in  $T$ are isotopic  if and only 
if their corresponding mixed  links in $S^3$  are, through an  isotopy
that  keeps the  specified unknotted  component  pointwise fixed.   By
applying to an oriented mixed link an  appropriate braiding we can then
turn it   into a `mixed'   braid  (a braid  that   keeps the specified
component pointwise fixed in the first position),  so that the closure
of this  braid is isotopic  to our mixed  link. An example of  a mixed
braid  is illustrated in the Introduction.  The set of  all  mixed braids on $n$
strands (where  the numbering excludes the first   fixed one) form the
group $B_{1,n}$, the   geometric version of  $\tilde{W}_n$.  Moreover,
similarly to  braid equivalence    in   $S^3$ we also    have   Markov
equivalence for mixed braids. (For details and proofs of the above the
reader is  referred to  \cite{L1} or \cite{L2}.)  Namely we  have  the
following.

\addtocounter{thm}{1}
\begin{thm} {\rm (cf.\ \cite{L2}, Theorem~3.)}\\
Let $L_1$,  $L_2$  be  two oriented  links  in $T$   and  ${\beta}_1$,
${\beta}_2$  be  mixed   braids  in  $\bigcup_{n=1}^{\infty}  B_{1,n}$
corresponding to $L_1$, $L_2$. Then $L_1$ is  isotopic to $L_2$ in $T$
if and   only  if    ${\beta}_1$  is  equivalent  to   ${\beta}_2$  in
$\bigcup_{n=1}^{\infty} B_{1,n}$  under  equivalence generated by  the
braid relations together with the following two moves:
\begin{itemize}
\item[(i)] Conjugation: If $\alpha,\beta \in B_{1,n}$
then $\alpha \sim \beta ^{-1}\alpha\beta$.
\item[(ii)] Markov moves: If $\alpha \in B_{1,n}$ then $\alpha \sim
\alpha{s_n}^{\pm 1}\in B_{1,n+1}$. 
\end{itemize}
\end{thm}

\noindent As already noted in  Section~1, there is  a strong resemblance  
between the Markov  moves and the special  property of a Markov trace.
Let  now $\pi$  denote    the canonical  quotient  map   $A\tilde{W}_n
\rightarrow H_n$ given in (3.1), and denote the generators of $H_n$ by
$t,g_1,\ldots,g_{n-1}$  as above. Let  also $\tau:H=\bigcup_{n \geq 0}
H_n \rightarrow A$ be the Markov trace (with parameter~$z \in A$) with
initial conditions $\tau(t_0't_1't_2' \cdots t_{k-1}')=y_k$ for all $k
\geq 1$.  Then   a braid  in $B_{1,n}$  can  be mapped  through  $\tau
\circ \pi$ to  an expression in  the variables $q, Q,  z, y_1, y_2, \ldots$.  
For an element $\alpha \in  B_{1,n}$ we shall denote by $\widehat{\alpha}$ 
its  closure. Then,  according to Theorem   5.2, in order to obtain an 
isotopy invariant $X$ for oriented  knots in $T$ we only need to  normalize  
$\tau$  so that  

\[X (\widehat{\alpha})=X(\widehat{\alpha{s_n}})=X(\widehat{\alpha s_n^{-1}}).\] 

This normalization has been done in \cite{L2}, (5.1), where   Jones's
normalization of Ocneanu's trace (cf.  \cite{Jones}) was followed. For
this purpose, we have to take some care  in the choice  of $A$ and the
parameter $z$.  We let  $A$ be  the field  of rational functions  over
$\QQ$  in   indeterminates $\sqrt{\lambda},\sqrt{q},\sqrt{Q},y_1,y_2,
\ldots$, and we let \[ z:=\frac{1-q}{q\lambda-1}.\]
(The reason for having square roots of $q$ and $Q$ will become clear
in the recursive formulae in (5.4) below; a square root of $\lambda$
is already required in the normalization of $\tau$.)

\begin{defn}{\rm
(cf.\  \cite{L2}, Definition 1.)}  \rm  For $\alpha$, $\tau$, $\pi$ as
above    let  \[X_{\widehat{\alpha}}  =    X_{\widehat{\alpha}}(q,  Q,
\sqrt{\lambda},  y_1, y_2,   \ldots)  =  
\Bigl[-\frac{1-\lambda q}{\sqrt{\lambda}(1-q)} \Bigr]^{n-1}
(\sqrt{\lambda})^e\, \tau(\pi(\alpha)), \] 
where $e$ is the exponent sum of the $s_i$'s that appear
in  $\alpha$.  (As noted in \cite{L2},  the $t'_i$'s can be ignored in
the   estimation   of   $e$    as  they  do  not    affect   it.) Then
$X_{\hat{\alpha}}$   depends    only  on     the   isotopy  class   of
$\widehat{\alpha}$, as a mixed  link representing an oriented  link in
$T$.
\end{defn}
If we look at $y_1,y_2, \ldots$ as parameters then this Definition  supplies
a family of invariants and Theorem~4.3 implies  that these are the only 
possible analogues of the 2-variable Jones polynomial for  oriented knots 
inside a  solid torus, which are related to the Iwahori-Hecke algebras of 
type $B$.  Note also that, if $\alpha \in B_{1,n}$ is a product of  the 
generators $s_i$ or their inverses (i.e, $\alpha$ does not involve the 
generator $t$) then $X_{\hat{\alpha}}$ is a rational function  of 
$\sqrt{\lambda}$ and $q$ only, and  it is exactly the   same as the  
invariant in \cite{Jones}, Definition~6.1. Geometrically, this means  that 
if an oriented knot in $T$ can  be enclosed in a 3-ball then the above 
invariant applied to this knot  will yield the  2-variable Jones polynomial 
(homfly-pt) for the knot, seen as a knot in $S^3$.

\medskip
5.4.  We shall now show  how to interprete  the above in terms of knot
diagrams,   and how to calculate  alternatively  the above solid torus
knot invariants using {\it initial conditions} and applying {\it skein
relations} on the mixed link diagrams.

Let $L_+ ,   L_-  , L_0$ be   oriented  mixed link diagrams that   are
identical, except in one crossing, where they are as depicted below:

$$\vbox{\picill3inby0.8in(mein4)  }$$

\noindent Let also $M_+ , M_- , M_0$  be oriented mixed link diagrams  
that are identical, except in the regions depicted below:

$$\vbox{\picill3inby0.8in(mein5)  }$$

\noindent In \cite{L2}, (5.2) it is shown that the knot invariant defined 
there (which is a  special member of  the family of invariants defined
above), satisfies the two recursive linear formulae.
\[\frac{1}{\sqrt{q}\sqrt{\lambda}}\,
X_{L_+} - \sqrt{q}\sqrt{\lambda}\,
 X_{L_-}  = (\sqrt{q}-\frac{1}{\sqrt{q}})\,X_{L_0} \ \ \ \  \]
\[\frac{1}{\sqrt{Q}}\, X_{M_+} - {\sqrt{Q}}\, X_{M_-}  =
(\sqrt{Q}-\frac{1}{\sqrt{Q}})\,X_{M_0} \ \ \ \ \]
These are the two skein relations that derive from the defining 
quadratic  relations of    $H$, and  the first   one  of  them is  the
well-known skein rule used   for   the evaluation of   the   homfly-pt
polynomial. The same reasoning applies to any invariant of  Definition 5.3 
above.

Take   now  a braid $\alpha$  in  $B_{1,n}$.  Then $\pi(\alpha)$ is an
$A$-linear combination of elements in the basis of $H_n$. (In terms of
diagrams, we have used  on $X_{\widehat{\alpha}}$ the skein  relations
above.)  This   fact   and Proposition   3.4 imply   now   that, under
conjugation  and the skein   relations, $X_{\widehat{\alpha}}$  can be
further  written as an $A$-linear combination  of the values of $X$ on
diagrams of  the form    $\widehat{d_1^{*}  \cdots d_n^{*}}$,    where
$d_i^{*}$ is  either $1, s_i$ or  $t'_{i,0}$ (with inverses mixed up).
If now $b$ is the number of $t'_{i,0}$'s  in $d_1^{*} \cdots d_n^{*}$,
then $\widehat{d_1^{*} \cdots  d_n^{*}}$   is  a  mixed link  on   $b$
components, such that each component links once  (in a positive sense)
with the special fixed one, but they are  otherwise unlinked with each
other. Therefore $\widehat{d_1^{*} \cdots d_n^{*}}$ is isotopic to the
mixed link $\widehat{t' t'_1 \cdots t'_{b-1}}$, pictured below.

$$\vbox{\picill1.2inby1.1in(mein6)  }$$

\noindent Now $\tau(t't'_1 \cdots t'_{k-1})=y_{k}$ is one of the initial
conditions from Theorem 4.3, and using Definition~5.3 we can
calculate $X_{\widehat{t' t'_1  \cdots t'_{k-1}}}$. Hence, we proved that 
the two skein rules together with the following (infinitely  many) initial 
conditions \[ X_{\widehat{1}}=1 \;\;\;\;(1 \in B_{1,1}),\;\;\; \quad  
\;\;\; X_{\widehat{\alpha_k}}=\Bigl[-\frac{1-\lambda q}{\sqrt{\lambda}(1-q)}
\Bigr]^{k-1}y_k\] (with $\alpha_k=t't_1' \cdots t_{k-1}' \in B_{1,k}$ for 
all $k \geq 1$) determine uniquely the invariant~$X$.

\addtocounter{thm}{1}
\begin{rems}{\rm  
Let $\alpha \in B_{1,n}$. The above discussion shows {\em geometrically}
that, firstly,  $\tau \circ \pi(\alpha)$ can be calculated  as a linear 
combination of terms $\tau(d_1^{*} \cdots d_n^{*})$ with $d_i^{*} 
\in {\cal D}_{n,0}'$ and, secondly, that $\tau(d_1^{*} \cdots d_n^{*})=
z^ay_b$ where $a$ is the number of $s_i$'s and $b$ is the number of
conjugates of $t$ in this element. This is the exact counterpart of the 
purely algebraic argument given before in (4.4).

Also, notice that the set of mixed links of the form $\widehat{t' t'_1
\cdots t'_{k-1}}$ form the basis of the submodule  of the 3rd skein module
of the solid torus (as calculated by Turaev in \cite{T}  and by Hoste and  
Kidwell in \cite{HK}), that is related to the Iwahori-Hecke algebras of 
type $B$.

Finally, the defining equation in Definition~5.3 already shows    that
$X_{\widehat{\alpha}}$ is a polynomial in $Q^{\pm},y_1,y_2,\ldots$. If
we  also  perform the   change  of variables  $x:=\sqrt{q\lambda}$ and
$r:=\sqrt{q}-\frac{1}{\sqrt{q}}$  then the   first skein   rule can be
rewritten  as  \[\frac{1}{x}  X_{L_+}-x   X_{L_-}= rX_{L_0}.\]   As in
\cite{Jones}, Proposition~6.2,   this  allows   us  to   deduce   that
$X_{\widehat{\alpha}}$  also is a Laurent  polynomial in the variables
$x$ and $r$.} \end{rems}

\begin{exmp} {\rm We shall now evaluate explicitly $X_{\widehat{\alpha}}$  
for some special choices of $\alpha$ as a Laurent polynomial in 
$x,r,Q,y_1,y_2,\ldots$.

\medskip
\noindent (a) If $\alpha$ is a product of generators $s_i$ or their inverses
(and does not involve~$t$) then $X_{\widehat{\alpha}}$ equals the known
2-variable invariant for oriented knots inside $S^3$ (see the remarks 
following Definition~5.3).

\medskip
\noindent (b) Consider $\alpha=s_1s_2^{-1}s_1ts_1^{-1}s_2$ and $\alpha'=ts_2$ 
in $B_{1,3}$ (see Remarks~2.6(b)). Using (4.4) we directly find that
$\tau \circ \pi (\alpha)=\tau \circ \pi(\alpha')=zy_1$. The exponent sum
$e$ of the factors $s_i$ equals $1$ in both cases. Therefore, we also
have \begin{eqnarray*}  X_{\widehat{\alpha}}=X_{\widehat{\alpha}'} & = & 
\left(\frac{1-\lambda q}{\sqrt{\lambda}(1-q)} \right)^2 \sqrt{\lambda}zy_1 \\
& = & \frac{\lambda q-1}{\sqrt{\lambda}(1-q)} \sqrt{\lambda}y_1 
\;\;\;\;\;\;\;\;\;\;\;\;\; \mbox{(inserting $z$)}\\
& = & \frac{1-x^2}{xr}y_1 
\;\;\;\;\;\;\;\;\;\;\;\;\;\;\;\;\;\;\;\;\;\;\;\; \mbox{(inserting 
$x$ and $r$).}
\end{eqnarray*}

\medskip
\noindent (c) For a general braid $\alpha \in B_{1,n}$, we first have to
consider its image in $H_n$ and express it as a linear combination in the 
standard basis of $H_n$:
\[h:=\pi(\alpha)=\sum_{w \in W} a_w g_w \quad \mbox{ with $a_w \in A$.}\]
Then we could compute the set $I(h)$ defined in Proposition~3.4 and 
finally use the recipe given in (4.4) to evaluate the trace $\tau(h)$.
Now the computation of the set $I(h)$ involves performing a base change
from the standard basis of $H_n$ to the new basis consisting of the elements
$r_1' \cdots r_n'$, with $r_i' \in {\cal R}_i'$ (see (3.2) and the proof of
Proposition~3.4). In practice, however, this will be quite cumbersome. 
A more economic way is by using the class polynomials of \cite{GP}. 

Recall from [{\em loc.\ cit.}] that for each $w \in W_n$ 
there exist elements $f_{w,C} \in {\ZZ}[q,Q]$ such that $\varphi(g_w)=
\sum_C f_{w,C} \varphi(g_{w_C})$ (sum over all conjugacy classes $C$ of
$W_n$), for all trace functions $\varphi$ on $H_n$. In [{\em loc.\ cit.}], 
Section~1, there is also given a recursive formula for computing $f_{w,C}$.  
Assume then that we know (for our given braid $\alpha$) the coefficients 
$a_w \in A$ and the class polynomials $f_{w,C}$ for all $w$ such that 
$a_w \neq 0$. Then we have
\[ \tau \circ \pi(\alpha)=\sum_C \left( \sum_{w \in W_n} a_wf_{w,C} 
\right) \tau(g_{w_C}).\]
Thus, we are reduced to calculating, once and for all, the values of
$\tau$ on basis elements corresponding to representatives $w_C$ 
of minimal length in the conjugacy classes $C$ of $W_n$. These classes
are parametrized by pairs of partitions $(\pi_1,\pi_2)$ such that the total sum
of the parts of $\pi_1$ and $\pi_2$ equals $n$ and where each part of $\pi_1$
(respectively $\pi_2$)  corresponds to a negative (respectively positive)
block. The formulae become more complicated as the number of negative blocks
involved in $w_C$ gets larger. Below we give the trace values for $n \leq 4$
and all $w_C$ which contain at most three negative blocks. (The only class 
$C$ for which we don't give the value is the one with representative 
$tt_1t_2t_3$; for each $n \geq 2$, we consider only those $w_C$ which are 
not already contained in $B_{1,n-1}$.) 
\[\begin{array}{cl}
B_{1,1}: & \tau(1)=1,\;\;\; \tau(t)=y_1.\\ & \\
B_{1,2}: & \tau(g_1)=z, \;\;\;\tau(tg_1)=zy_1,\\
& \tau(tt_1)=((q-1)(Q-1)y_1+(q-1)Q)z+qy_2.\\ & \\
B_{1,3}: & \tau(tg_2)=zy_1,\;\;\; \tau(g_1g_2)=z^2,\;\;\; 
\tau(tg_1g_2)=z^2y_1,\\
& \tau(tt_1g_2)=z\tau(tt_1) \;\;\;\;\;\; \mbox{ (see $B_{1,2}$)}, \\
& \tau(tt_1t_2)=((q-1)(Q^2-Q+1)y_1+Q(q-1)(Q-1))(q^3-1)z^2\\
& \;\;\;\;\;\;\;\;\;\;\;\;\;\;\;\; +(Qy_1+(Q-1)y_2)(q^3-1)qz+q^3y_3. 
\end{array}\]
\[\begin{array}{cl}
B_{1,4}: & \tau(g_1g_3)=z^2,\;\;\;\tau(g_1g_2g_3)=z^3,\\
& \tau(tg_1g_3)=\tau(tg_2g_3)=z^2y_1,\;\;\;
\tau(tg_1g_2g_3)=z^3y_1, \\
& \tau(tt_1g_3)=z\tau(tt_1),\;\;\;\tau(tt_1g_2g_3)=z^2\tau(tt_1) 
\;\;\;\;\;\; \mbox{ (see $B_{1,2}$)},\\
& \tau(tg_1t_2g_3)=((Q-1)y_1+Q)(q-1)(q^2+1)z^3+q(q^2-q+1)y_2z^2,\\
& \tau(tt_1t_2g_3)=z\tau(tt_1t_2) \;\;\;\;\;\; \mbox{ (see $B_{1,3}$)},\\
& \tau(tt_1t_2t_3) \;\;\;\;\; \mbox{ (four negative blocks)}. \end{array}\]
Let us now consider the example of the braid $\alpha \in B_{1,5}$ given in
 the Introduction. By first using the defining properties for $\tau$, we can
eliminate the generators $g_4$ and $g_3$ and are reduced to computing \[\tau(\alpha)
=z^2\tau(g_2^{-1}g_1^{-3}tg_1g_2^{-1}g_1g_2^{-1}t^{-2}g_1^{-1}).\]
In principle, this can be done by following the above scheme (but the
result will be quite cumbersome and we do not want to print it here). For a 
similar computation see \cite{L2}, Example pp. 237. Based on the programs 
in \cite{GP} it is straightforward to implement the above algorithmic 
description in a computer program.} \end{exmp}

\bigbreak \noindent {\bf Acknowledgements.} The second author would like 
to acknowledge financial support by the European Union and to thank 
Professor tom Dieck for an invitation to the SFB 170 in G\"{o}ttingen, 
where part of this work was completed. 


\bigbreak
\noindent {M.G.: UFR de Math\'ematiques, Universit\'e Paris 7, 2 Place Jussieu,
75251 Paris, France, E-mail geck@mathp7.jussieu.fr}

\noindent {S.L.: Mathematisches Institut, G\"ottingen Universit\"at,
 3--5 Bunsenstrasse, 37073 G\"ottingen, Germany, E-mail
sofia@cfgauss.uni-math.gwdg.de}

\end{document}